\documentclass[11pt]{amsart}
\usepackage{amsmath,amsthm,amsfonts,amscd,amssymb}
\usepackage[mathscr]{eucal}
\usepackage[all]{xy}
\usepackage[colorlinks]{hyperref}

\newtheorem{thm}{Theorem}
\newtheorem{lemma}[thm]{Lemma}

\newtheorem{proposition}[thm]{Proposition}
\newtheorem{corollary}[thm]{Corollary}

\newtheorem{conjecture}[thm]{Conjecture}

\newcommand{\p}{{\mathfrak p}}

\renewcommand{\P}{{\mathbb P}}

\newcommand{\F}{{\mathcal F}}

\newcommand{\J}{{\mathscr J}}

\newcommand{\Q}{{\mathbb Q}}

\newcommand{\x}{{\times}}

\newcommand{\la}{{\langle}}
\newcommand{\ra}{{\rangle}}
\newcommand{\lra}{{\,\longrightarrow\,}}

\renewcommand{\th}{{\text{th}}}

\newcommand{\lla}{{\, \longleftarrow \,}}

\DeclareMathOperator{\chara}{char}

\DeclareMathOperator{\Ass}{Ass}

\title{Minimal Primes Over Permanental Ideals}
\author{George A. Kirkup}
\address{University of California, Berkeley}
\date{September 16, 2005}
\email{kirkup@math.berkeley.edu}

\begin{document}

\begin{abstract}
In this paper we discuss minimal primes over permanental ideals
of generic matrices.  We give a complete list of the minimal primes 
over ideals of $3 \x 3$ permanents of a generic matrix, and show that
there are monomials in the ideal of maximal permanents of a $d \x 2d-1$ 
matrix if the characteristic of the ground field is sufficiently large.
We also discuss the Alon-Jaeger-Tarsi Conjecture, using our results
and techniques to strenthen the previously known results.
\end{abstract}

\maketitle

\section{Introduction}
\subsection{Definitions}

Let $A = (a_{i,j})$ be any $n \x n$ matrix over a commutative ring $R$.  
The permanent of this matrix is
$$\sum_{\sigma \in S_n} a_{1,\sigma(1)}a_{2,\sigma(2)} \cdots a_{n,\sigma(n)}.$$
This is often described as the formula for the determinant of 
$A$ without the minus signs.  

Fix $m,n,d$ and a field $k$.  We define $I_d(m,n) \subset k[x_{i,j}]$ to be the 
ideal generated by the permanents of the $d \x d$ submatrices of 
$M_{m,n} = (x_{i,j})$, a generic $m \x n$ matrix over $k$.  The main goal 
of this paper is to understand the minimal primes over $I_d(m,n)$.

\subsection{Background}

Two recent papers discuss this notion.  In \cite{LS}, Laubenbacher and Swanson
carefully analyze $I_2(m,n)$ for all $m,n$.  They use the language of ideals, 
minimal primes, radical, Gr\"obner bases, \emph{etc.} to give a very complete 
understanding of these ideals.

In \cite{Yu}, Yu discusses some relationships between the rank of a matrix and 
the ``perrank'' of a matrix, namely the size of the largest submatrix whose 
permanent is nonzero.  He does not use the language of ideals, preferring an 
entirely set-theoretic approach to begin to answer some combinatorial questions.
These were raised in \cite{AT} Alon and Tarsi in 1989 and again treated in 
\cite{BBLS} in 1994.

\subsection{Overview}

Our techniques will be very algebraic, and will be simple applications of the
multilinearity of permanents which is discussed in Section \ref{multi}.  In
Section \ref{2sec} we review relevant results from \cite{LS} and give several
proofs based on the results of the previous section.

In Section \ref{monsec} we discuss the main conjecture of Chapter, which is
\begin{conjecture}\label{preconj}
If $\chara k > d$ or $\chara k = 0$ and $n \geq 2d-1$ then the minimal 
primes containing $I_d(m,n)$ must either contain a column of the generic 
matrix or the $d-1 \x d-1$ permanents of some $m-1$ rows.  
\end{conjecture}
This would inductively give the minimal primes over the ideals $I_d(m,n)$ for 
large $m,n$ if we knew the minimal primes over $I_d(m,n)$ for $m,n \leq 2d-1$.
We prove Conjecture \ref{preconj} for the case $d=3$ in Section \ref{35sec}.  
What we are able to show in general is that any prime over $I_d(d,2d-1)$ which
contains no entry from some row must contain the $(d-1)\x (d-1)$ permanents 
of the other $d-1$ rows.  

We continue in Section \ref{d+1sec} where we discuss the ideals\linebreak 
$I_d(d,d+1)$ in detail.  Then in Section \ref{34sec} we focus on $I_3(3,4)$, paying 
special attention to $3 \x 4$ matrices with no entries vanishing whose maximal permanents 
do vanish.  

In Section \ref{submaxsec} we discuss in general the case in which $m,n > d$.  We apply
these ideas to the case $d=3$ in Section \ref{44sec}, in which we list the
minimal primes over $I_3(4,4)$.  From this we deduce the minimal primes over
$I_3(m,n)$ for all $m,n$.

In Section \ref{AJTsec} we relate our conjecture and results to the combinatorial 
conjectures in \cite{AT}.

In this section we strengthen the results of \cite{AT}
and \cite{BBLS}.  

\section{The Multilinearity of the Permanent and the Algebra of Permanents}\label{multi}

Let $A = (a_{i,j})$ be a $d \x d$ matrix.  Then the permanent of $M$
can be expressed as
$$\sum_{1}^{d} x_{i,1} \cdot A^{\hat {\imath}}_{\hat 1}$$
where $A^{\hat {\imath}}_{\hat 1}$ is the $(d-1)\x(d-1)$ permanent of the 
submatrix of $A$ obtained by omitting row $i$ and column $1$.
This is similar to the expansion by minors (Laplace expansion) for
determinants.  This expansion can obviously done for any row or column
and this multilinearity is the key observation for our entire study
of permanental ideals.

We begin with an example of the importance of multilinearity that we
will use throughout this Chapter.
If $A$ is a $d \x (d-1)$ matrix, and column vectors
$v_1, \dotsc, v_d$ have the property that $(A \medspace v_i)$ has
permanent $0$ for all $i$, then the vector $(A^{\hat 1}, \dotsc, A^{\hat d})$,
of the $(d-1)\x(d-1)$ permanents of $A$, is in the kernel of the
matrix $(v_1 \medspace v_2 \dots v_d)$.  Thus, by rank-nullity, either the
determinant of $(v_1 \medspace v_2 \dots v_d)$ is $0$ or the
$(d-1)\x(d-1)$ permanents of $A$ vanish identically.

Another application concerns the ``algebra of permanents.''
Fix $d < n$ and consider the ring $S  = k[A_{i_1,\dotsc,i_d}]$ where 
$1 \leq i_j \leq n$ for each $j$ and the $A_{i_1,\dotsc,i_d}$ are
indeterminants.  Let $R = k[x_{r,s}]$ where
$1 \leq r \leq d$ and $1\leq s \leq n$.  Then there is a map
$\phi \colon S \lra R$ taking $A_{i_1,\dotsc,i_d}$ to the permanent of the
matrix $(c_{i_1}, \medspace \dotsc, \medspace c_{i_d})$ where
$c_j$ is the $j^{\th}$ column of the generic $d \x n$ matrix $(x_{r,s})$.
The image of the map $S \lra R$ is the algebra of maximal permanents,
and is isomorphic to $S/(\ker \phi)$.  We can find interesting elements 
of $\ker \phi$ in the following way.

Let $H$ be a matrix whose columns are indexed by multi-subsets, $\J$, of 
$\{1,\dotsc,n\}$ whose cardinality is $d-1$.  Then let the column of
$H$ indexed by $\J$ be the vector $A_{1,\J}, \dotsc, A_{d,\J}$.
For example, if $d=3$, and $n=4$ then
$$H = \begin{pmatrix} A_{1,1,1}&A_{1,1,2}&A_{1,1,3}&A_{1,1,4}&A_{1,2,2}&A_{1,2,3}&A_{1,2,4}&A_{1,3,3}&A_{1,3,4}&A_{1,4,4} \\
                      A_{1,1,2}&A_{1,2,2}&A_{1,2,3}&A_{1,2,4}&A_{2,2,2}&A_{2,2,3}&A_{2,2,4}&A_{2,3,3}&A_{2,3,4}&A_{2,4,4} \\
                      A_{1,1,3}&A_{1,2,3}&A_{1,3,3}&A_{1,3,4}&A_{2,2,3}&A_{2,3,3}&A_{2,3,4}&A_{3,3,3}&A_{3,3,4}&A_{3,4,4} \\
                      A_{1,1,4}&A_{1,2,4}&A_{1,3,4}&A_{1,4,4}&A_{2,2,4}&A_{2,3,4}&A_{2,4,4}&A_{3,3,4}&A_{3,4,4}&A_{4,4,4} \end{pmatrix}.$$
We now focus our attention on a particular column, with index $\J$.
Then let $A^{\hat {\imath}}_{\J}$ be the permanent of the $(d-1)\x(d-1)$ 
matrix with columns $\J$ and omitting row $i$.  Then the $\J^{\th}$ 
column can be expanded as
$$\begin{pmatrix} \sum_{i} x_{i,1} A^{\hat {\imath}}_{\J}\\
                  \sum_{i} x_{i,2} A^{\hat {\imath}}_{\J}\\
                  \vdots                          \\
                  \sum_{i} x_{i,n} A^{\hat {\imath}}_{\J}\end{pmatrix}.$$
Therefore, no matter what $\J$ is, the column above is in the span of
$$\begin{pmatrix} x_{1,1}\\
                  x_{1,2}\\
                  \vdots \\
                  x_{1,n}\end{pmatrix}, \dotsc, \begin{pmatrix} x_{d,1}\\ x_{d,2}\\ \vdots  \\ x_{d,n}\end{pmatrix}$$
so the rank of $H$ is at most $d$.  Therefore, the $(d+1)\x (d+1)$
minors of $H$ vanish.

\section{$I_2(2,3)$}\label{2sec}
This section reviews a result of Laubenbacher and Swanson \cite{LS}, and
gives three different proofs of it, suggesting three different
approaches to the problem in general.
\begin{thm}\label{23thm}
If $\chara k \neq 2$ then any prime containing $I_2(2,3)$
either contains an entire row of $M_{2,3}$ or it contains 
an entire column of $M_{2,3}$.
\end{thm}
\begin{proof}
Let $\p$ be a prime over $I_2(2,3)$.
The first argument we give is closely related to the computational
proof in \cite{LS}.  If there is no entry of the first row of $M_{2,3}$
in $\p$ then since $I_{2,3}$ is homogeneous in each row and each column, 
we can dehomogenize each column
so that the $x_{1,1}=x_{1,2}=x_{1,3}=1$.  Then $I_2(2,3)$ becomes
$\la x_{2,1}+x_{2,2},x_{2,1}+x_{2,3},x_{2,2}+x_{3,2}\ra$.  Since
$\chara k \neq 2$, this ideal is generated by 
$\la x_{2,1},x_{2,2},x_{2,3}\ra$.  Therefore, there must
be some entry of $M_{2,3}$ in $\p$.  We may reindex so that $x_{1,1}\in \p$.
Thus $\p$ must also contain $x_{2,1} x_{1,2}$ and
$x_{2,1} x_{1,3}$ so either $\p$ contains the first row
or the first column of $M_{2,3}$.

The next two proofs use our notation from the previous section,
$$A_{i,j} = x_{1,i}x_{2,j}+x_{1,j}x_{2,i}.$$
Since $A_{1,2} = A_{1,3} = 0$ modulo $I_2(2,3)$, the vector $(x_{1,1},x_{2,1})$
is in the kernel of the matrix of $1\x 1$ permanents
$$T_{2,3} = \begin{pmatrix} x_{1,2} & x_{1,3} \\ x_{2,2} & x_{2,3} \end{pmatrix}.$$
Therefore, $\la x_{1,1}, x_{2,1}\ra \cdot \la \det T_{2,3} \ra \in I_2(2,3)$
so $\p$ either contains a column of  $M_{2,3}$ or the $2\x 2$ minors of $M_{2,3}$.  
Since $\chara k \neq 2$, we can combine the minors with 
the permanents to form all polynomials of the form $x_{1,i}x_{2,j}$ with
$i \neq j$.  The result then follows similarly to the end of the first proof.

The last proof uses the fact that the symmetric matrix
$$H_{2,3} = \begin{pmatrix} A_{1,1} & A_{1,2} & A_{1,3}\\
                            A_{1,2} & A_{2,2} & A_{2,3}\\
                            A_{1,3} & A_{2,3} & A_{3,3}\end{pmatrix}$$
has rank $2$ by the discussion in Section \ref{multi}.  Since $A_{i,j} =0$ 
modulo $I_2(2,3)$ whenever $i\neq j$, the determinant of $H_{2,3}$ is 
$A_{1,1}\cdot A_{2,2}\cdot A_{3,3}$.  Therefore, one of the $A_{i,i}$
must be in $\p$.  Since $\chara k \neq 2$, $A_{i,i} \in \p$
implies that $x_{1,i}\cdot x_{2,i} \in \p$ and thus one of these entries is
in $\p$.  The proof again follows in the same manner as the first proof.
\end{proof}

These three proofs of the result reflect three different
perspectives on how to proceed with the study of permanental
ideals.

\section{Monomials in $I_d (d,2d-1)$}\label{monsec}
In this section we will focus on $d \x (2d-1)$ matrices whose 
maximal permanents vanish, showing that there must be some entries 
which vanish (under mild hypotheses) and giving some general tools 
for understanding the minimal primes over the ideal $I_d(d,2d-1)$.  
We use the approach in the first proof of Theorem \ref{23thm}.
The case of $d \x n$ matrices with $n < 2d-1$ is less understood.

The discussion in this section will revolve around the following 
conjecture, which will be proved in the case $d=3$ in Section
\ref{35sec}

\begin{conjecture}\label{mainconj}
Let $k$ be any field with $\chara k > d \text{ or } 0$ and $M_{d,2d-1}$ 
be the generic $d \x (2d-1)$ matrix over $k$.  Any minimal prime
over $I_d(d,2d-1)$ either contains a column of $M_{d,2d-1}$
or contains the $(d-1)\x (d-1)$ permanents of some set of $d-1$ 
rows of $M_{d,2d-1}.$  If $\chara k \leq d$ then, in addition to these
possibilities, a minimal prime may instead contain the $d\x d$ minors of 
$M_{d,2d-1}$.
\end{conjecture}

First we set up some notation for the section.  Let $M_{d,n} = (x_{i,j})$
be the generic $d \x n$ matrix where $n > d$.  If
$$\overline \alpha = (\alpha_1, \dotsc , \alpha_{d-1})$$
is a strictly increasing sequence with $1 \leq \alpha_1$ and
$\alpha_{d-1} \leq 2d-1$ then let $A_{\overline \alpha}$ be
the $(d-1) \x (d-1)$ subpermanent of $M$ obtained by omitting
the first row and including columns $\alpha_1, \dotsc , \alpha_{d-1}$.
Now define
$$A'_{\overline \alpha} = A_{\overline \alpha} \cdot \prod_{i \not\in \overline \alpha} x_{1,i}.$$

If 
$$\overline \beta = (\beta_1, \dotsc , \beta_{d})$$
is a strictly increasing sequence with $1 \leq \beta_1$ and
$\beta_{d} \leq n$ then let $P_{\overline \beta}$ be
the $d \x d$ subpermanent of $M$ obtained by including columns 
$\beta_1, \dotsc , \beta_{d}$, and as above, let
$$P'_{\overline \beta} = P_{\overline \beta} \cdot \prod_{i \not\in \overline \beta} x_{1,i}.$$
Note that
$$P'_{\overline \beta} = \sum_{\overline \alpha \subset \overline \beta} A'_{\overline \alpha}.$$
The next proposition will use sums of permanents heavily, so let
$$S_i(\overline\alpha) = \sum P'_{\overline \beta}$$
where the sum runs over all $\overline \beta$ such that the cardinality
of $\overline \beta \cap  \overline \alpha$ is $i$.  Finally, let
$$T_i(\overline\alpha) = \sum A'_{\overline \alpha'}$$
where the sum runs over all $\overline \beta$ such that the cardinality
of $\overline \alpha' \cap  \overline \alpha$ is $i$.

\begin{proposition}\label{preevid}
For any $d,n$ with $2d-1 > n > d$ and any $\overline \alpha$ of size $d-1$,
the ideal $I_d(d,n)$ contains
$$(n-d+1)!(A'_{\overline \alpha} + (-1)^{n-d} T_{2d-n-2}(\overline \alpha)).$$
\end{proposition}
\begin{proof}
In this proof we will treat $A'_{\overline \alpha}$ as atomic.  We will prove
the result by constructing an element of $I_d(d,n)$ and checking the number of
times each $A'_{\overline \alpha}$ appears.

Fix $\overline \alpha$ and let $S_i = S_i(\overline \alpha)$.
Let
\begin{equation}\label{incexc}
B = \sum_{i = 2d-n-1}^{d-1} (-1)^i (i+n-2d+1)! (d-i-1)! S_i.
\end{equation}
We claim that $B = \pm (n-d+1)!(A'_{\overline \alpha} + T_{2d-n-2})$.

Each $S_i$ is a sum of various $A'_{{\overline {\alpha}}'}$.  In fact,
if the cardinality of ${\overline {\alpha}}' \cap \overline \alpha$ is $i$,
then $S_i$ contains $n-(d-1)-(d-i-1)=n-2d+i+2$ copies of $A'_{{\overline {\alpha}}'}$.
On the other hand, $S_{i+1}$ also contains $d-1-i$ copies of
$A'_{{\overline {\alpha}}'}$.

Given any $A'_{{\overline {\alpha}}'}$, for which the cardinality of 
${\overline {\alpha}}' \cap \overline \alpha$ is $i$, it can only appear in the
definition of $S_i$ or $S_{i+1}$.  Thus we can measure its contribution
to $B$ from the definition in \eqref{incexc}.  It appears
$$(-1)^{i+1} (n-2d+i+2)! (d-i-2)! (d-1-i) + (-1)^i (n-2d+i+1)! (d-1-i)! (n-2d+i+2) = 0$$
times.
Therefore, the only contribution to $B$ from $A'_{{\overline {\alpha}}'}$
comes from the ${\overline {\alpha}}'$ such that the cardinality of 
$\overline \alpha \cap \overline \alpha'$ is $2d-n-2$ or $\alpha' = \alpha$.  In
the former case, $A'_{{\overline {\alpha}}'}$ appears in $S_{2d-n-1}$ 
$$d-1-(2d-n-2) = n-d+1$$
times.  Therefore, the total number of times $A'_{{\overline {\alpha}}'}$ 
appears in the definition of $B$ is 
\begin{equation*}
\begin{split}
(-1)^{2d-n-1} (2d-n-1+n-2d+1)! & \cdot (d-(2d-n-1)-1)! \cdot (n-d+1) \\
                               & = (-1)^{n+1} (0)!\cdot (n-d)! \cdot (n-d+1) \\
                               & = (-1)^{n+1}(n-d+1)!.
\end{split}
\end{equation*}
$A'_{\overline\alpha}$ appears in $S_{d-1}$ $n-(d-1)$ times 
and in no other $S_i$.  Therefore, the number of times $A'_{\overline\alpha}$
appears in the definition of $B$ is
\begin{equation*}
\begin{split}
(-1)^{d-1} (d-1+n-2d+1)! & \cdot (d-(d-1)-1)! \cdot (n-(d-1)) \\
                         &= (-1)^{d-1} (n-d)! \cdot (0)! \cdot (n-d+1) \\
                         &= (-1)^{d-1} (n-d+1)!.
\end{split}
\end{equation*}
This completes the proof since $B$ is a sum of polynomials
in $I_d(d,n)$
\end{proof}

The following proposition is proved in exactly the same manner, and
shows why the ideal $I_d(d,2d-1)$ is of special importance.

\begin{proposition}\label{evid}
$d! \cdot A'_{\overline \alpha} \in I_d(d,2d-1)$ for any $d, \overline \alpha$.
\end{proposition}
\begin{proof}
We can use the same proof as in Proposition \ref{preevid}.
Let
\begin{equation}\label{incexc1}
B = \sum_{i = 0}^{d-1} (-1)^i i! (d-i-1)! S_i.
\end{equation}
We claim that $B = d! \cdot A'_{\overline \alpha}$.

Given any $A'_{{\overline {\alpha}}'}$, with the cardinality of 
${\overline {\alpha}}' \cap \overline \alpha$ $i$, it only appears
in the definition of $S_i$ and $S_{i+1}$.  Thus we can count the number
of times it appears in $B$ from the definition \eqref{incexc1}; it appears
$$(-1)^{i+1} (i+1)! (d-i-2)! (d-1-i) + (-1)^i i! (d-1-i)! (i+1) = 0$$
times.  Therefore, the only contribution to $B$ from an $A'_{{\overline {\alpha}}'}$
comes from $A'_{\overline {\alpha}}$, since there is no $i = d$.  Thus
$B = d! \cdot A'_{\overline \alpha}$ as claimed.
\end{proof}

Notice that this generalizes the first proof of Theorem \ref{23thm}.
We get the following corollaries.

\begin{corollary}\label{moncor}
If $\chara k > d$ then $\prod x_{i,j} \in I_d(d,2d-1)$.
\end{corollary}
\begin{proof}
We induct on $d$.  Certainly $x_{1,1} \in I_1(1,1)$.

Now assume that $\prod x_{i,j} \in I_{d-1}(d-1,2d-3)$.
Then it is clear that $\prod x_{i,j} \in I_{d-1}(d-1,2d-1)$.
However, the proposition implies that for any 
$f \in I_{d-1}(d-1,2d-1)$, 
$$f \cdot \prod x_{d,j} \in I_{d,2d-1},$$
proving the result.
\end{proof}

The following corollary states Proposition \ref{evid} set-theoretically.
\begin{corollary}\label{setevid}
Suppose that $\chara k > d$ and $M$ is a $d \x (2d-1)$ matrix over $k$
whose maximal subpermanents vanish.  Then for any row, either 
there is an entry in that row which vanishes, or the $(d-1) \x (d-1)$
subpermanents vanish in the other $d-1$ rows.
\end{corollary}

\section{The Ideal $I_3(3,5)$}\label{35sec}

In this section, we will use the results of the previous section
to prove Conjecture \ref{mainconj} in the case $d=3$.  If 
$\chara k = 2$, $I_3(3,5)$ is just the determinantal ideal, which is prime.
Therefore, for the rest of this section we assume that $\chara k \neq 2$.

For the sake of clarity, we review the results from the previous section
in the case $d = 3$.  Let $A_{i,j}$ be the permanent of the $2 \x 2$
submatrix of $M_{3,5}$, with columns $i,j$ and rows $2,3$.  Then
Proposition \ref{evid} implies that if $\chara k = 0$ or $\chara k > 3$
then
$$A_{1,2} x_{1,3} x_{1,4} x_{1,5} \in I_3(3,5).$$
Moreover, Proposition \ref{preevid} implies that for all $k$
(recall $\chara k \neq 2$),
$$A_{1,2} x_{1,3} x_{1,4} - A_{3,4} x_{1,1}x_{1,2} \in I_3(3,5).$$
All polynomials of these form which can be obtained by permuting
the rows and columns of $M$ are also in $I_3(3,5)$.

\begin{thm}
If $d =3$ and $\chara k \neq 2,3$ then any minimal prime containing 
$I_3(3,5)$ either contains a column of $M$ or the $2\x 2$ permanents 
of some two rows.  If $\chara k = 3$ any prime not containing one of
these must contain the $2\x2$ minors of some collection of $4$ columns.
\end{thm}
\begin{proof}
Let $\p$ be a prime over $I_3(3,5)$.  For any three columns of $M$, 
either $\p$ contains the determinant of those three columns, or it 
contains the $2 \x 2$ permanents of the other $2$ columns by the 
multilinearity of the permanents.  Suppose that there is a $3 \x 3$ 
submatrix whose determinant is not in $\p$, and reindex so they are 
columns $3,4,5$. Then the $2\x2$ permanents of the first $2$ columns 
are in $\p$.  Therefore, by Theorem \ref{23thm} either a row or a 
column of this $3 \x 2$ matrix is in $\p$.  If a column is in $\p$ we 
are done.  If a row is in $\p$ then we can reindex so $x_{1,1},x_{1,2} \in \p$.  
If $x_{1,i}\in \p$
for some other $i$, reindex so it $i=3$, and then the $2\x 2$
permanents $A_{1,2},A_{1,3},A_{2,3}$ (omitting row $1$) are in $\p$
or $x_{1,4},x_{1,5} \in \p$.  In the latter case, we are done, and in
former case, we are also done by Theorem \ref{23thm}.  Therefore, we can
assume that $x_{1,i} \not\in \p$ for $i=3,4,5$.  In this case, by
Proposition \ref{preevid}, we again have $A_{1,2},A_{1,3},A_{2,3} \in \p$,
which proves the claim.

There remains the case in which the $3\x 3$ minors of $M_{3,5}$ are in $\p$.  
In this case, the matrix $H_{3,5}$ from Section \ref{multi} also must have rank 
$2$ (modulo $\p$), so its $3\x 3$ minors are in $\p$.  Consider three 
well-selected columns of $H$,
$$H' = \begin{pmatrix} A_{1,1,2} & A_{1,1,3} & 0         \\
                       A_{1,2,2} & 0         & A_{2,2,3} \\
                       0         & A_{1,3,3} & A_{2,2,3} \\
                       0         & 0         & 0         \\
                       0         & 0         & 0         \end{pmatrix}.$$
The existence of $H'$ shows that if the $2 \x 2$ minors
of the last two rows are not in $\p$ then all the entries of $H'$ are in 
$\p$.  This implies that $3\x 3$ permanents of the matrix
\begin{equation}\label{AJTmatrix3}
\begin{pmatrix} x_{1,1} & x_{1,2} & x_{1,3} & x_{1,1} & x_{1,2} & x_{1,3} \\
                x_{2,1} & x_{2,2} & x_{2,3} & x_{2,1} & x_{2,2} & x_{2,3} \\
                x_{3,1} & x_{3,2} & x_{3,3} & x_{3,1} & x_{3,2} & x_{3,3} \end{pmatrix}
\end{equation}
vanish.  If $\chara k \neq 2,3$, then $3! \neq 0$, 
so we may assume that there are at least two entries
from the first row in $\p$ by Proposition \ref{evid}.  Thus, by the previous paragraph,
if $\chara k \neq 3$ we have proved the result.

If $\chara k =3$, we can still use Proposition \ref{preevid}.  Let
$A_{i,j}$ be the $2 \x 2$ subpermanent of $M_{3,5}$ obtained by omitting
the first row and including columns $i,j$ (possibly $i=j$).  We are in the
situation in which the maximal minors of $M_{3,5}$ are in $\p$, but
the $2\x 2$ minors of the last two rows are not in $\p$, which implied that
maximal permanents of matrix \eqref{AJTmatrix3} were in $\p$.  If an entry
of \eqref{AJTmatrix3} were in $\p$ then the previous paragraph would imply the
result.  Therefore, we can dehomogenize each column
(set $x_{1,1}=x_{1,2}=x_{1,3}=1$) so Proposition \ref{preevid} implies
that $A_{i,j}=A_{k,l}$ for all $i,j,k,l \leq 3$ by applying the proposition
to each set of $4$ columns of \eqref{AJTmatrix3}.  Therefore, $A_{i,j} - A_{i,k}=0$
for all $i,j,k$ so the permanent of
$$\begin{pmatrix} x_{2,i} & x_{2,j} - x_{2,k} \\
                  x_{3,i} & x_{3,j} - x_{3,k} \end{pmatrix}$$
vanishes for all choices of $i,j,k$.  Therefore, either $x_{2,j} = x_{2,k}$ for all 
$j,k$ or the $2\x 2$ minors of the bottom two rows of \ref{AJTmatrix3}
vanish.  Since this is true for all $j,k$, the $2\x 2$ minors of the 
bottom two rows vanish in either case.  Since $I_3(3,5)$ is also homogeneous
in each row and none of the entries of \eqref{AJTmatrix3} vanish,we can dehomogenize 
the bottom two rows so that $x_{2,1} = x_{3,1} = 1$.  Therefore we have specialized
the matrix \eqref{AJTmatrix3} to
\begin{equation}
\begin{pmatrix} 1 & 1   & 1   & 1 & 1   & 1 \\
                1 & x_2 & x_3 & 1 & x_2 & x_3 \\
                1 & x_2 & x_3 & 1 & x_2 & x_3 \end{pmatrix}.
\end{equation}
The $3\x 3$ permanents of this matrix yield the equations
$$2(1 + x_2 + x_2) = 2(1 + x_3 + x_3) = 0$$
so $x_2 = x_3 = 1$.  Rehomogenizing, this means that
the $2 \x 2$ minors of matrix \eqref{AJTmatrix3} are in $\p$.

Therefore, for any pair of columns, either the $2\x 2$ minors
of those two columns are in $\p$ or the $2 \x 2$ minors of the
other three columns are in $\p$.  This implies that either
there is either a column vanishes modulo $\p$ or the $2 \x 2$
minors of four columns of $M_{3,5}$ are in $\p$.  This completes
the proof.
\end{proof}

\section{Minimal Primes over $I_d(d,d+1)$}\label{d+1sec}

Corollary \ref{moncor} implies that if $\chara k > d$ and if the 
maximal permanents of a $d \x (2d-1)$ matrix vanish, then at least
one entry of the matrix is $0$.  The natural question is whether
this is sharp.  In particular, if the permanents of a
$d \x n$ matrix, with $n<2d-1$ vanish, does there need to be an
entry of the matrix which is $0$?  The following example shows that
at least in the case $d \x (d+1)$, there need not be any entries
which vanish.  These examples were constructed by making many of the
entries in the matrices $1$ and solving the resulting equations.
$$\begin{pmatrix} 1      & \cdots & 1      & 1      & 2-3d       \\
                  \vdots & \ddots & \vdots & \vdots & \vdots     \\
                  1      & \cdots & 1      & 1      & 2-3d       \\
                  1      & \cdots & 1      & 2-2d   & (d-2)(d-1) \\
                  1      & \cdots & 1      & d      & (2d-1)(d-2)  \end{pmatrix}$$
For example, if $d =3$, this is
$$\begin{pmatrix} 1 & 1 & 1  & -7  \\
                  1 & 1 & -4 & 2   \\
                  1 & 1 & 3  & 5   \end{pmatrix}.$$

Moreover, the following $4 \x 6$ matrix (where $i^2 = -1)$ has its maximal 
permanents vanish,
$$\begin{pmatrix} 1 & 1  & 1 & 1  & 1  & 1   \\
                  1 & 1  & 1 & -1 & -1 & -1  \\
                  1 & 1  & 1 & i  & i  & i   \\
                  1 & 1  & 1 & -i & -i & -i  \end{pmatrix}.$$
At this point I am not sure whether there are any $4 \x 6$ matrices
over $k$ with vanishing maximal permanents and no vanishing entries
if $k$ does not have a square root of $-1$.

We now discuss $I_d(d,d+1)$, the smallest example of an ideals defined
by $d \x d$ permanents, in detail.  The next section will explicitly
go through the calculations in Section \ref{multi} in this case.

Let $M = (x_{i,j})$ be a $d \x (d+1)$ generic matrix and let $P_j$ be 
the permanent of the square matrix obtained by omitting the
$j^{\th}$ column of $M$.  Let
$$B_{r,p,q} = \frac{\partial}{\partial x_{r,p}} P_q.$$
Notice that $B_{r,p,q} = B_{r,q,p}$ and $B_{r,p,p} = 0$.
Moreover, for $p \neq q$, $B_{r,p,q}$ is the $d-1 \x d-1$ 
subpermanent of $M$ obtained by omitting row $r$ and columns $p,q$.

Now consider the matrix of partial derivatives
$$L'_p = \begin{pmatrix} B_{1,1,p}   & \cdots & B_{1,d+1,p}   \\
                         \vdots      &        & \vdots      \\
                         B_{d,d+1,p} & \cdots & B_{d,d+1,p} \end{pmatrix}$$
and let $L_p$ be the $d \x d$ matrix obtained by omitting
the $p^{\th}$ column of $L'_p$, which is $0$.  We define $f_p$ to 
be the determinant of $L_p$.

Another natural construction for a matrix of partial derivatives
holds a row $r$ constant:
$$W_r = \begin{pmatrix} B_{r,1,1}   & \cdots  & B_{r,1,d+1} \\
                  \vdots      &         & \vdots \\
                  B_{r,d+1,1} & \cdots  & B_{r,d+1,d+1} \end{pmatrix}.$$
Let $g_r$ be the determinant of the symmetric matrix $W_r$.

Notice the similarity between $W_r$ and $H_{2,4}$ from Section \ref{multi}.
However, there is an important difference.  Namely,  $A_{i,j}$ in
$H_{2,4}$ is the permanent of the $2 \x 2$ matrix which has columns
$i,j$ whereas $B_{r,i,j}$ \emph{omits} columns $i,j$.

\begin{proposition}\label{structJ}
For any $i,j$, we have $x_{i,j} f_j , x_{i,j} g_i \in I_d(d,d+1)$.
\end{proposition}
\begin{proof}

By symmetry we may assume that $i = j = 1$.

For $x_{1,1}f_1$, note that for all $j \neq 1$,
$$P_j = \sum_i x_{i,1} B_{i,1,j}.$$
Let $e_j$ be the determinant of the $d-1 \x d-1$ submatrix
of $L_1$, omitting the first column and the $j^{\th}$ 
row.  Now consider the polynomial
\begin{equation*}
\begin{split}
f'_j &= \sum_j (-1)^j e_j P_j \\
     &= \sum_j (-1)^j e_j \sum_i x_{i,1} B_{i,1,j} \\
     &= \sum_i x_{i,1} \sum_j (-1)^j B_{i,1,j} e_j.
\end{split}
\end{equation*}
For $i = 1$, the interior sum is just the expansion by minors
of $\det L_1 = f_1$.  However, for $i \neq 1$, the interior
sum is the expansion by minors of the determinant of $L_1$
where the first column is replaced by the $i^{\th}$.  Thus,
for $i \neq 1$, the interior sum is $0$.
Thus $f'_j = x_{1,1} f_1$, so $x_{1,1}f_1 \in I_d(d,d+1)$.

The second statement is proved in exactly the same method,
noting that
$$P_j = \sum_p x_{1,p} B_{1,p,j}.$$
(Recall that if $j = p$ then $B_{1,p,j} = 0$.)  From here,
the proof proceeds exactly as the previous part.
\end{proof}

The previous proposition allows us to make some statements about
the minimal primes over $I_d(d,d+1)$.  Let 
$$J_d = I_d(d,d+1) + \la f_j, g_i \ra \subset k[x_{r,s}].$$

\begin{corollary}\label{mindd+1}
If $Q$ is any minimal prime containing $I_d(d,d+1)$, then either
$Q$ contains a row, a column, or contains $J_d$.
\end{corollary}
\begin{proof}
This is immediate from the previous proposition.  If $Q$ does not
contain some $f_j$, then it contains $x_{i,j}$ for all $i$.  Thus
$Q$ contains the prime ideal generated by $x_{i,j}$ for all $i$
and the $d \x d$ permanent $P_j$.

On the other hand, if $Q$ does not contain $g_i$ for some $i$, then
$Q$ contains the ideal generated by $x_{i,j}$ for all $j$.  So $Q$
contains a row.
\end{proof}

\section{$J_3$, A Prime Containing $I_3 (3,4)$}\label{34sec}
In this section, we assume that the characteristic of $k$ is $0$, although I believe
that the results hold (except for the rational point discussion) for all odd
characteristics.

For ease of notation, let
$$M = \begin{pmatrix} x_1 & x_2 & x_3 & x_4 \\
                      y_1 & y_2 & y_3 & y_4 \\
                      z_1 & z_2 & z_3 & z_4 \end{pmatrix}$$
and let $I_3 (3,4)$ be the ideal of $3\times 3$ subpermanents of $M$.

\begin{proposition}\label{34prop}
Let $k = \Q$.  The ideal $J_3$, defined in Section \ref{d+1sec} is a prime 
ideal of codimension $4$, and $I_3(3,4)$ is the intersection of $J_3$, 
the primes containing a row of $M$, and those containing a column of $M$.  
Therefore, $I_3(3,4)$ is a radical complete intersection.
\end{proposition}

The proofs are due to Singular \cite{GPS01} and Macaulay 2 \cite{M2}.  Singular 
will compute the primary decomposition of $I_3 (3,4)$ (using the GTZ algorithm -- 
the SY algorithm did not terminate) in several hours.

The most interesting ideal in $\Ass(I_3 (3,4))$ is certainly $J_3$.
Since $I_3(3,4)$ is a complete intersection, we can compute its degree as 
$$\deg(J_3) = \deg(I_3(3,4)) - 3\deg(P_1) - 4\deg(P_3) = 81 - 3(1) - 4(3) = 66$$
where $P_1$ is a prime containing a row of $M$  and $P_3$ is a prime 
containing a column of $M$, which can be verified directly by Macaulay 2 \cite{M2}.
The singular locus of $I_3(3,4)$ is contained in the variety defined by $J_3$. 

In the previous section we exhibited a rational point on the variety defined by $J_3$.
The question of finding all rational points with all entries non-zero (full support) can 
be simplified greatly by the observation that multiplication of any row or column by a scalar 
results in the multiplication of the permanent by that scalar.  This is a result of the
homogeneity of $I_d(m,n)$.  In particular, we may  normalize five entries to obtain a 
matrix of the form
$$M = \begin{pmatrix} 1 & 1 & 1 & 1 \\
                      u & a & b & c \\
                      u & d & e & f \end{pmatrix}.$$

This is a Zariski open set in the variety defined by $J_3$, so there are
few matrices which cannot be written in this form, and 
most of those are contained in other irreducible components of $I_3(3,4)$.  Let 
$J'$ be the ideal generated by the $3 \times 3$ subpermanents of $M$, which is 
still homogeneous.  The variety defined by $J'$ is certainly Zariski dense in that 
defined by $J_3$ since it contains the open set of all matrices with full support.

One interesting situation arises when one of the entries, say $a = -u$.  This is evidently
a special case since $a = -u$ implies that the $2 \x 2$ subpermanent in the top left corner
is $0$.  It turns out that the ideal $I_3(M) + \la a + u \ra$ has four primary components, 
and the irreducible varieties they correspond to may be represented in the following ways.
$$\begin{pmatrix} 1 & 1 & 1 & 1    \\
                  0 & 0 & b & c    \\
                  0 & 0 & eb & -ce \end{pmatrix}, 
\begin{pmatrix}   1 & 1 & 1 & 1  \\
                  0 & 0 & 0 & 0  \\
                  0 & d & e & f  \end{pmatrix},$$
$$\begin{pmatrix} 1 & 1  & 1 & 1 \\
                  u & -u & 0 & 0 \\
                  u & u  & e & -e \end{pmatrix}$$
and
$$\begin{pmatrix} 1 & 1  & 1 & 1 \\
                  u & -u & u & u \\
                  u & 0  & e & -u-e \end{pmatrix}$$
Notice that in no case can the matrix have full support.

Thus, we can change coordinates, $A = a+u,B = b+u,\ldots F = f+u$ and still invert the 
variables.  The four $3 \times 3$ subpermanents of $M$ can be thought of as relations
between the $2 \times 2$ subpermanents of the matrix consisting of the last two rows of $M$.
In particular, after some reorganizing, we can rewrite the generators of $J'$ as
$$J' = \la AE+BD -2u^2, AF+CD-2u^2,BF+CE-2u^2,u\cdot(A+B+C+D+E+F-6u) \ra.$$

Since we can invert $A, B, C$, we get the relations
$$E = \frac{2u^2 - BD}{A} \;\;\;\;\;\;\; F = \frac{2 - CD}{A}$$
and thus
$$D = \frac{u^2\cdot(-A + B + C)}{BC}.$$
By symmetry, this means that given the $u,A,B,C$ coordinates, the rest of the 
coordinates are determined (by field operations), as long as $A,B,C,u$ are nonzero.
Recall, however, that one of the generators of $J_3$, $g_3$, is a relation only among
the elements of the first two rows of $M$.  In particular, we have the following
relation
\begin{equation*}
\begin{split}
h &= u^2ab+u^2ac+u^2bc+ua^2b+ua^2c+uab^2+uac^2+ub^2c+ubc^2+a^2bc+ab^2c+abc^2 \\
  &= u^2(ab+ac+bc)+u(a^2b+a^2c+ab^2+ac^2+b^2c+bc^2)+(a^2bc+ab^2c+abc^2) \\
  &= \sigma_1 \sigma_3 - \sigma_4
\end{split}
\end{equation*}
where $\sigma_1 = a+b+c+u$, $\sigma_3 = abc+abu+acu+bcu$, and $\sigma_4 = abcu$.

The form $h$ defines a quartic surface in $\P^3$ which is symmetric in all four variables and
has fourteen double points.  Since it is quadratic in each variable, we can view it as a
double cover of $\P^2$.  In that light, we can find several rational
curves on the surface, by adding constraints which make the discriminant vanish.  
For example, if $abc(a+b+c) = 0$, then the constant term vanishes, when viewed as 
a quadratic in $u$.

We can also view the surface as an elliptic surface over $\P^1$, since $h$ is a cubic over
the function field of $\P^1$ (pick any $2$ variables to distinguish).  In this light, if
we can find a rational point on that elliptic curve (corresponding to a rational curve),
we can multiply it by integers using the group law on the elliptic curve.  This will 
yield infinitely many rational curves on the surface unless the original point was a torsion
point.

Notice that if $A = B$ then the constant term is $0$, so there is a unique solution for $C$
in terms of $A$ and $B$.  In particular, we note that the following matrix has its $3 \x 3$
subpermanents vanish for any choice of $a$.
$$\begin{pmatrix} 1   & 1  & 1  & a+2      \\
                  1   & a  & a  & -a(2a+1) \\
                  a+1 & -a & -a & a(1-a)   \end{pmatrix}$$
Moreover, for all $a \neq -2,-1,-1/2,0,1$, this yields a matrix with full support that has
its $3 \x 3$ subpermanents vanish.

Another interesting special case is when $A = 1/B$ since then the leading 
term vanishes. In that case we can write another one-parameter family of matrices 
in $J$,
$$\begin{pmatrix} 1   & 1        & a+1    & 3(a+1)           \\
                  1   & a        & -a     & a^2+a+1          \\
                  a+2 & -a(2a+1) & a(1-a) & (a-1)(2a+1)(a+2) \end{pmatrix}.$$
In this case, we get a matrix of full support unless $a = -2,-1,-1/2,0,1$ or $a$
is a root of $x^2+x+1$ or $3=0$.

Other information about $J_3$ that seems interesting is derived from the free resolution, 
which Macaulay 2 calculates as
$$0 \lla S/J_3 \lla S \lla S^{11} \lla S^{34} \lla S^{42} \lla S^{24} \lla S^{6} \lla 0$$
with Betti diagram
$$\begin{matrix}
        \text{total}&  :&   1&   11&   34&   42&   24&   6\\
                   0&  :&   1&    .&    .&    .&    .&   .\\
                   1&  :&   .&    .&    .&    .&    .&   .\\
                   2&  :&   .&    4&    .&    .&    .&   .\\
                   3&  :&   .&    .&    .&    .&    .&   .\\
                   4&  :&   .&    .&    6&    .&    .&   .\\
                   5&  :&   .&    4&   12&   12&    .&   .\\
                   6&  :&   .&    .&    .&    .&    .&   .\\
                   7&  :&   .&    3&   16&   30&   24&   6\end{matrix}$$

The Castelnuovo-Mumford regularity of $J_3$ is $7$ and
its Hilbert polynomial is
$$\frac{11}{840} z^7 + \frac{1}{15} z^6 + \frac{1}{30} z^5 + \frac{41}{12} z^4 + \frac{653}{120} z^3 + \frac{1171}{60} z^2 - \frac{2647}{210} z + 7.$$

\section{Submaximal Permanents}\label{submaxsec}
In this section we will introduce some notation that will make the discussion
in Section \ref{44sec} clearer.  

Fix $d,m,n$ and let $v$ be a weakly decreasing sequence of positive 
integers of finite length $l(v)$, for example $(3,2,2)$.  Let $s(v)$ 
be the sum of all the entries of $v$ and define 
$$d'(v) = s(v) - l(v) +1$$
and $e(v) = d - d'(v)$.

Suppose that $s(v) \leq m,n$ and $d'(v) \leq d$.  Then
we define a prime of type $v$ as follows.  Let $M'$ be any submatrix of $M$ of 
size $m' \x n'$ where $m' + n' + e = m +n$.  Then consider any ordered partition of
the rows and columns of $M'$ such that the sizes of the sets in the partition (both
for rows and columns) matches the entries of $v$, except perhaps for extra one
element sets in the partition.  Let $a_i$ be the permanent of the $v_i \x v_i$ 
submatrix of $M'$ given by the rows and columns in the sets of size $v_i$ in the 
partitions.  Then let $\p$ be generated by $a_i$ for all $i$ and all the entries of $M'$
not used in any of the $a_i$.  Any prime that can be obtained in this way is called
a prime of type $v$.

Notice that given a sequence $v$, if $v' = (v,1)$ then the primes of type $v$ are
the same as the primes of type $v'$ (so long as $s(v') \leq m,n$) so we can somewhat
simplify the description of the types by assuming that no entry of $v$ is $1$, except
in the case of the empty sequence, which we notate by $(1)$.

For some easy examples of this, suppose that $d = 2$.  Primes of type $(1)$ are
those which are generated by all the elements of $M$ except for one row or one
column.  Primes of type $(2)$ are those which are generated by some $2 \times 2$
subpermanent and all the entries of $M$ outside this $2 \times 2$ block.  There
are no other kinds of primes of type $v$ in the case $d = 2$ since $d'(v)$ can
be at most $2$.

\begin{proposition}
Fix $m,n,d$, and suppose that $v$ is a weakly decreasing sequence of positive
integers of finite length.  If $s(v) \leq m,n$ and $d'(v) \leq d$ then primes
of type $v$ contain $I_d(m,n)$.
\end{proposition}
\begin{proof}
Since it is clear that the primes of type $v$ are indeed prime, we only need to
show that the variety of such a prime is contained in the algebraic set defined
by $I_d(m,n)$.  The next lemma, which is a corollary of the expansion by minors
of permanents allows us to reduce to the case in which $d'(v) = d$.

\begin{lemma}\label{exminlem}
Suppose that $m_1 \leq m$, $n_1 \leq n$, and $d_1 \leq d$.  If 
$$d - d_1 \geq (m+n) - (m_1 +n_1)$$
then $I_{d_1}(m_1,n_1) \supset I_d(m,n)$.
\end{lemma}

Thus we assume $d'(v) = d$, and we induct on $l(v)$.  If $l(v) = 1$, then
$d = v_1$ and a prime of type $v$ is generated by a $d \x d$ subpermanent
of $M$ and all the other entries of $M$, which evidently contains $I_d(m,n)$.

Let $v = (v_1 \ldots v_l)$, let $\overline v = (v_1 \ldots v_{l-1})$, and 
let $\p$ be any prime of type $v$.  We can rearrange columns and rows so that
$\p$ is generated by the upper-left $v_l \x v_l$ subpermanent, all the
other entries from the first $v_l$ rows and columns, and in the lower right
$m-v_l \x n -v_l$ submatrix, $\overline M$, is a prime of type $\overline v$.  Let 
$$\overline d = d'(\overline v) = d - v_l -1.$$
By induction, $\p$ contains all the $\overline d \x \overline d$ subpermanents
of $\overline M$.  Now, consider any $d \x d$ submatrix, $T$, of $M$.  If $T$
contains the first $v_l$ rows of $M$, then since $\p$ contains all the
$v_l \x v_l$ subpermanents of the submatrix given by the first $v_l$ rows
of $T$, expansion by minors implies that the permanent of $T$ is
contained in $\p$.  Now assume that $T$ does not contain the first $v_l$
rows of $M$.  Then it must contain at least $d - v_l + 1$ rows from the
rest of $M$.  However, by induction, $\p$ contains the $\overline d \x \overline d$ 
subpermanents of this lower submatrix.  Now we can use Lemma \ref{exminlem}
to finish the proof.  In this case $m_1 \geq d - v_l -1$, $m = n = n_1 = d$, and
$d_1 = \overline d= d - v_l - 1$.  Therefore, since
$$d - d_1 = d - (d - v_l - 1) \geq (m + n) - (m_1 + n_1)$$
$\p$ contains the permanent of $T$.  Thus we have shown that $\p$ contains
$I_d(m,n)$.
\end{proof}

\section{The Ideals $I_3(m,n)$}\label{44sec}
We have already completed our discussion of the minimal primes over $I_3(3,n)$.
In this section we will discuss the minimal primes over $I_3(4,4)$.  The primes
from the previous section which are relevant are those of type (1), (2), (3),
and (2,2).  In addition, we say a prime is of type (3A) if it contains the 
ideal $J_3$ (from section \ref{34sec}) for some $3 \x 4$ (resp. $4 \x 3$) 
submatrix and all the entries of the remaining rows (resp. columns).

\begin{thm}\label{44thm}
If $\chara k > 3$ then the minimal primes over $I_3(4,4)$ are of type 
(1), (2), (3), (3A), or (2,2).  If $\chara k =3$ then the ideal of $2 \x 2$
minors of the generic matrix is the only other minimal prime.
\end{thm}
\begin{proof}
If $\p$ is any prime over $I_3(4,4) = I$ and $\p$ contains all the
entries from some row or column, then $\p$ obviously contains a prime
of type (3) or (3A).  If $\p$ contains three elements from some row
or column (say $x_{1,1},x_{1,2},x_{1,3}$) then $\p$ either contains
$x_{1,4}$ or the $2\x 2$ permanents of the matrix
$$\begin{pmatrix} x_{2,1} & x_{2,2} & x_{2,3} \\
                  x_{3,1} & x_{3,2} & x_{3,3} \\
                  x_{4,1} & x_{4,2} & x_{4,3} \end{pmatrix}$$
and thus contains a prime of type (2) or (1).

If $\p$ contains two elements from some row or column (say $x_{1,1},x_{1,2}$),
then $\p$ either contains $x_{1,3},x_{1,4}$ or the $2 \x 2$ permanents of the
first two columns.  By Theorem \ref{23thm} this means $\p$ either contains
a column of $M_{4,4}$ or the entries in a $2\x 2$ block of $M_{4,4}$.  The
former case has been dealt with already.  In the latter case, suppose
$\p$ contains $x_{1,1},x_{1,2},x_{2,1},x_{2,2}$.  Let the $A_{i,j}$ be the 
$2 \x 2$ permanent of the bottom two rows and columns $i,j$.  We have the
following relations from the $3\x 3$ permanents,
$$x_{1,3} A_{1,4} + x_{1,4} A_{1,3}, x_{2,3} A_{1,4} + x_{2,4} A_{1,3}$$
so either $x_{1,3}x_{2,4} - x_{1,4}x_{2,3} \in \p$ or $A_{1,3},A_{1,4} \in \p$.
However, since $x_{1,3}x_{2,4} - x_{1,4}x_{2,3} \in \p$ already, the former
case would imply that there were three elements from some row in $\p$ (since
$\chara k \neq 2$).  In the latter case, we already have $A_{1,2} \in \p$, and
similarly $A_{2,3},A_{2,4} \in \p$.  Therefore, either $\p$ contains a prime of
type (2,2) or $\p$ contains three entries from one of the bottom two rows.

If $\p$ contains only one element from some row (say $x_{1,1}$) then by Proposition
\ref{preevid} the $2\x 2$ permanents of the first $2$ columns and last three rows
vanish.  Therefore, we can use one of the previous cases.

Finally, if $\p$ contains no entries of $M$, then we can dehomogenize each column so
$x_{1,i} = 1$ for each $i$.  As above let $A_{i,j}$ be the $2 \x 2$ permanents in the
last two columns.  By Proposition \ref{preevid} we know $A_{1,2} = A_{3,4}$,
$A_{1,3} = A_{2,4}$, and $A_{1,4} = A_{2,3}$.  If the $2\x 2$ minors of the first two
rows do not vanish, then by the multilinearity of permanents, the $3\x 3$ minors of
$$\begin{pmatrix} 0       & A_{3,4} & A_{2,4} & A_{2,3} \\
                  A_{3,4} & 0       & A_{1,4} & A_{1,3} \\
                  A_{2,4} & A_{1,4} & 0       & A_{1,2} \\
                  A_{2,3} & A_{1,3} & A_{1,2} & 0       \end{pmatrix}.$$
The top left $3\x 3$ minor of this matrix is $2 A_{1,2}A_{1,3}A_{1,4}$.
Therefore, one of those is in $\p$, say $A_{1,2}$.  This implies that
$A_{3,4} \in \p$ as well.  Moreover, this means that $A_{1,3}+A_{2,3}=0$
so substituting in the matrix $H$ from Section \ref{multi} the the $3 \x 3$
minors of the following matrix vanish.
$$\begin{pmatrix} A_{1,1}  & 0        & A_{1,3}  & -A_{1,3} \\
                  0        & A_{2,2}  & -A_{1,3} & A_{1,3}  \\
                  A_{1,3}  & -A_{1,3} & A_{3,3}  & 0        \\
                  -A_{1,3} & A_{1,3}  & 0        & A_{4,4}  \end{pmatrix}$$
Taking the first three columns and rows 1,2,4 the minor is
$$-A_{1,1} A_{1,3} (-A_{1,3}) - (-A_{1,3})A_{2,2}A_{1,3} = A_{1,3}^2(A_{1,1} + A_{2,2}).$$
Therefore, either $A_{1,3} = 0$ which would mean that the $2\x 2$ permanents
of the bottom two rows vanish, and hence an entry vanishes, or $A_{i,i} =0$ for
all $i$ which also implies that entries of the matrix vanish.  Since no entries of
the matrix vanish, this is a contradiction.  Therefore, the $2 \x 2$ minors of
every pair of rows vanishes, which is to say that the $2 \x 2$ minors of $M$ are
in $\p$.  If $\chara k \neq 3$, this implies that many entries of $M$ vanish, so
the result is proved.
\end{proof}

\begin{corollary}\label{Perm3Cor}
If $m,n \geq 4$ and $\chara k \neq 2,3$ then all minimal primes over $I_3(m,n)$ are 
of type (1), (2), (3), (3A), or (2,2).  If $\chara k =3$ then another minimal prime
is the ideal of $2\x 2$ minors of $M_{m,n}$.
\end{corollary}

Notice that if $m,n \geq 4$ then all the primes of type 1, 2, 3, 3A, and (2,2) are 
minimal over $I_3(m,n)$.

\section{The Alon--Jaeger--Tarsi Conjecture}\label{AJTsec}
Now we turn our attention to one of the major motivations
for the study of permanental ideals, the Alon--Jaeger--Tarsi
Conjecture.  We will briefly review the results of \cite{AT}
and then show how our results and conjectures apply to this
conjecture.

The Alon--Jaeger--Tarsi conjecture was stated in \cite{AT} as
\begin{conjecture}\label{AJTconj1}
Let $A$ be a nonsingular $d \x d$ matrix over a finite field
$k$ with cardinality $q \geq 4$.  There exists a vector $v$
in $k^n$ such that both $v$ and $Av$ have no zero component.
\end{conjecture}

In \cite{AT}, Alon and Tarsi show that if $A$ is a matrix over
$k$ such that when the matrix $A$ is repeated $q-2$ times, the
resulting matrix $A' = (A \mid A \mid \dotsb \mid A)$ has a
nonvanishing maximal permanent then there is such a vector.
They prove Conjecture \ref{AJTconj1} for $q = p^k$ where $p$ is prime and
$k \geq 2$.  In \cite{BBLS} the authors extend this result to
the cases $q \geq d+1 \geq 4$ and $q \geq d+2 \geq 3$ using a
simple counting argument.  In Corollary \ref{BBLSgen} we will
show this result for $q \geq d \geq 4$.

Let $A = (a_{i,j})$ be a $d \x n$ matrix over a field $k$.
The key objects of \cite{AT} are the polynomials
$$P_{A}(X_1,\dotsc,X_n) = \prod_{i=1}^d (\sum_{i=1}^n a_{i,j} X_j)$$
and
$$P'_{A}(X_1,\dotsc,X_n) = X_1 X_2 \cdots X_n \cdot P_A(X_1,\dotsc,X_d)$$
The definition of $P_A$ immediately implies
\begin{lemma}\label{polylink}
If $X = (X_1,\dotsc,X_n)$ is any vector in $k^n$ then $P'_A(X) = 0$ 
if and only if $X$ and $A(X)$ both have no nonzero entries.
\end{lemma}

We introduce the notation of \cite{AT} here for convenience and
consistency.  Let
$\alpha = (\alpha_1,\dotsc,\alpha_n)$ be a vector of nonnegative
entries whose sum is $d$.  Then let $c_{\alpha}$ be the coefficient
of the monomial $X_1^{\alpha_1}\dotsb X_n^{\alpha_n}$ and let
$A_{\alpha}$ be the $d \x d$ matrix with $\alpha_i$
repetitions of the $i^{\th}$ column of  $A$.  Claim 1 in
\cite{AT} is
\begin{lemma}\label{permlink}
The permanent of $A_{\alpha}$ is $c_{\alpha} \cdot \prod (\alpha_i)!$.
\end{lemma}

If one could show that for an invertible matrix $A$, the $d \x 2d$
matrix $A' = (A \mid A)$ has a non-vanishing maximal permanent, then
Conjecture \ref{AJTconj1} would follow.  Therefore, Jeff Kahn
conjectured that, indeed, $A'$ has a non-vanishing maximal permanent.
Yu discusses this conjecture and has some interesting theorems
related to it in \cite{Yu}. The following conjecture is a refinement 
of Kahn's by De Loera \cite{SPE}.  It should be noted that Conjecture 
\ref{mainconj} obviously implies this result set-theoretically.

\begin{conjecture}\label{AJTconj}
Let $M_{d,d}$ be the generic $d\x d$ matrix over a field $k$.
Then if $I$ is the ideal of $d \x d$ permanents of the $d \x 2d$
matrix $(M_{d,d} \mid M_{d,d})$, $(\det M_{d,d})^d \in I$.
\end{conjecture}

If $k$ has characteristic $2$, the permanent is the determinant
so the conjecture is obvious.  If the characteristic of the
field is $3$ then we have the following refinement, the proof of 
which is related to the proof of the result in \cite{AT}.

\begin{proposition}\label{char3AJT}
If $\chara k = 3$ and $I_d$ is the ideal of $d \x d$ permanents 
of the $d \x 2d$ matrix $(M_{d,d} \mid M_{d,d})$, then
$$\prod_{1}^d D_i(M_i) \subset I_d$$
where $M_i$ is the $i \x d$ submatrix of $M_{d,d}$ given by the
first $i$ rows and $D_i(M_i)$ is the ideal of $i \x i$ minors of
$M_i$.
\end{proposition}
\begin{proof}
Let $\{a_1,\dotsc,a_{d-1}\}$ be any multiset of columns of
$M_{d,d} = M$.  Then modulo $I_d$, the permanent of the matrix
with columns $\{i,a_1,\dotsc,a_{d-1}\}$ is $0$ for any $i$ because
if there is no number repeated three times, then the permanent is
a generator of $I_d$ and if there is a number repeated three times,
the $3\x 3$ permanents of those three columns vanish since $6 = 0$.
Therefore, by expansion by ``minors,'' the permanent vanishes.
This implies that $D_d(M) \cdot P_{d-1}(M \mid M) \in I_d$ where
$P_{d-1}(M \mid M)$ is the ideal of $(d-1) \x (d-1)$ subpermanents
of $(M \mid M)$.  Thus, we can induct, using the case $d=1$ which
is obvious as the base case.
\end{proof}

This refinement is \emph{not} true for arbitrary fields.  In fact,
a quick calculation using Macaulay 2 \cite{M2} shows that if $k = \Q$
then 
$$D_3(M_3) \cdot D_2(M_2) \cdot D_1(M_1) \not\subset I_d(M_3 \mid M_3).$$
Computations using Macaulay 2 also suggest that in characteristic $3$,
$$\prod_{j=1}^{d} D_i(M_d) \subset I_d(M_d \mid M_d),$$
which is stronger than Proposition \ref{char3AJT}.

The next lemma, whose proof is contained in the proof of the Proposition
\ref{char3AJT}, may be of great use in proving Conjecture \ref{AJTconj}, 
especially if one could classify the minimal primes over $I_{d-1}(d,d)$.

\begin{lemma}\label{AJTsublemma}
Let $M_{d,d}$ be the generic $d\x d$ matrix over a field $k$.
Then if $I$ is the ideal of $d \x d$ permanents of the $d \x 2d$
matrix $(M_{d,d} \mid M_{d,d})$,
$$(\det M_{d,d})\cdot I_{d-1}(d,d) \in I.$$
\end{lemma}

Using this Lemma and Corollary \ref{Perm3Cor}, we prove
\begin{thm}\label{BBLSgen}
Conjecture \ref{AJTconj1} is true for $q \geq d \geq 4$.
\end{thm}
\begin{proof}
Let $A$ be an invertible $d \x d$ matrix.  By Lemma
\ref{permlink}, there is a vector $X$ with only nonzero 
coordinates such that $AX$ has nonzero coordinates unless
the matrix $(A \mid \dotsb \mid A)$, with $A$ repeated
$q-2$ times has vanishing maximal permanents.  Therefore,
we can assume that this is the case.  Since $A$ is 
invertible, Lemma \ref{AJTsublemma} implies that the 
$d-1 \x d-1$ permanents of $(A \mid \dotsb \mid A)$, with
$A$ repeated $q-3$ times, vanish.  By repeated application
of Lemma \ref{AJTsublemma}, we see that the $d-q+3$ sized
permanents of $A$ vanish.  Since $q \geq d$, we know that
the $3 \x 3$ permanents of $A$ vanish.  However, by
Corollary \ref{Perm3Cor}, the only invertible $d \x d$
matrices with its $3 \x 3$ permanents vanishing are of
type (2,2) with $d=4$.  Therefore, we need only consider
$A$ of the form
$$\begin{pmatrix} x_1 & x_2 & 0   & 0   \\
                  y_1 & y_2 & 0   & 0   \\
                  0   & 0   & z_1 & z_2 \\
                  0   & 0   & w_1 & w_2 \end{pmatrix}$$
with the added constraint that $(A \mid A)$ has its $4 \x 4$
permanents vanishing.  This is impossible unless one of the
rows vanishes, which contradicts the assumption that $A$ is
invertible.
\end{proof}

Now we turn our attention away from invertible matrices, and
consider arbitrary matrices which have vectors with only
nonzero coordinates such that their images have the same property.

Lemma \ref{permlink} is useful whenever $\alpha_i <p$ for all
$i$.  Otherwise, $(\alpha_i)! = 0$ for some $i$ and thus the lemma 
is trivial for that $\alpha$.  To get past this, we would need to replace
the permanent of $A_{\alpha} = b_{i,j}$ by another object, defined similarly.
Fix $\alpha$ and let $S_d$ be the symmetric group on the set $\{1,\dotsc,d\}$.
Further, let $S_{\alpha}$ be the quotient group $S_d/(\prod S_{\alpha_i})$
where $S_{\alpha_1}$ acts on the set $\{1,\dotsc,\alpha_1\}$, $S_{\alpha_2}$
on the set $\{\alpha_1,\dotsc,\alpha_1+\alpha_2\}$, \emph{etc}.  Then the
product $b_{1,\sigma(1)} b_{2,\sigma(2)} \dotsb b_{d,\sigma(d)}$ is well-defined
in the sense that for any two representatives of $\sigma$ in $S_d$ the above 
products are equal.  Using this notation, we get the following lemma.
\begin{lemma}\label{betterlink}
$$c_{\alpha} = \sum_{\sigma \in S_d/(\prod S_{\alpha_i})} \prod_1^d b_{i,\sigma(i)}$$
\end{lemma}
This lemma is what was implicitly used in \cite{AT} to prove Lemma \ref{permlink}.

The last piece of \cite{AT} that we will use is
\begin{lemma}\label{reduceq}
Let $f$ be any polynomial in $k[x_1,\dotsc x_n]$ where $k$ is the 
finite field with $q$ elements, $k = \F_q$.  $f$ is identically 
$0$ over $k$ if and only if $f$ can be reduced to the zero polynomial
by the relations $x_i^q = x_i$.
\end{lemma}

We are now ready to prove
\begin{thm}\label{AJTlargechar}
Suppose that $A$ is any $d \x n$ matrix over the finite field 
$k = \F_q$ and $q > d+1$.  There is a vector $X \in k^n$ such
that neither $X$ nor $AX$ has any $0$ entries if and only if
no row of $A$ is identically $0$.
\end{thm}
\begin{proof}
It is clear that if a row of $A$ is identically $0$ then
every vector in the image of $A$ has an entry that is $0$.
Therefore, we may assume that there is no vectors $X$ such 
that neither $X$ nor $AX$ has any $0$ entries, and we will 
prove that $A$ has a row which is identically $0$.

By Lemma \ref{polylink} we know that the polynomial $P'_A$ 
is identically $0$ on $k^n$.  Since $q > d+1$, Lemma \ref{reduceq} 
implies that $c_{\alpha} = 0$ for every $\alpha$.

We proceed by induction on $d$, since the result is clear
for $1 \x n$ matrices.  For $\alpha = (d,0,\dotsc,0)$, 
$c_{\alpha} = 0$ implies that there is an entry of the first 
column of $A$ which is $0$.

We will prove that either a row of $A$ is identically $0$ or
the first column is identically $0$.  Assume that there are 
exactly $d>e>1$ entries of the first column which are $0$.  Reindex 
so that $a_{1,1} = \dotsb = a_{e,1} = 0$.  Then let 
$$L = \{ \alpha \mid \alpha_1 = d-e \}.$$
By Lemma \ref{betterlink}, for every $\alpha \in L$,
$$c_{\alpha} = a_{e+1,1} \dotsb a_{d,1} \cdot c'_{(\alpha_2,\dotsc,\alpha_n)}.$$
where $A'$ is the $e \x (n-1)$ submatrix of $A$, omitting
the first column and including the first $e$ rows and
$c'_{(\alpha_2,\dotsc,\alpha_n)}$ is the coefficient
of $X_2^{\alpha_2}\dotsc X_n^{\alpha_n}$ of $P'_{A'}$.
If $c'_{(\alpha_2,\dotsc,\alpha_n)} = 0$ for all $\alpha$
then by induction $A'$ must have a row that is identically
$0$, but that implies that $A$ also has such a row.  

Therefore, one of $a_{e+1,1} \dotsb a_{d,1}$
must be $0$, so more than $e$ entries of the first column vanish
which contradicts the assumption that exactly $e$ entries
vanish.  Therefore, the first column must be identically $0$, but
this argument implies that every column is identically $0$, so $A$
is identically $0$.  Therefore, in any case, $A$ has
a row which is identically $0$.
\end{proof}

Now we consider the case where $d+1$ is a prime, and use 
Corollary \ref{setevid} to establish the following theorem.

\begin{proposition}\label{AJTd+1}
Suppose that $A$ is any $d \x n$ matrix over the finite field 
$k = \F_{d+1}$ where $d+1$ is prime.  Further, suppose that 
there is no vector $X \in k^n$ such that neither $X$ nor $AX$ 
has any $0$ entries.  Then either a row of $A$ is identically $0$ 
or
$$A = \begin{pmatrix} 1      & b_1      & 0      & \dots & 0 \\
                      1      & b_2    & 0      & \dots & 0 \\
                      \vdots & \vdots & \vdots &       & \vdots \\
                      1      & b_d    & 0      & \dots & 0 \end{pmatrix}$$
up to scaling of the rows and permutation of the columns,
where $\{b_1,\dotsc,b_d\}$ is the multiplicative group $k^{*}$.
\end{proposition}
\begin{proof}
Since $P'_A$ is identically $0$, by Lemmas \ref{polylink} and
\ref{permlink}, the matrix $(A \mid \dotsb \mid A)$ where $A$
is repeated $d-1$ times has its maximal permanents vanish.  By
Corollary \ref{setevid} either all but $2d-2$ entries in every
row of this matrix vanish or the $(d-1)\x (d-1)$ permanents of
some $d-1$ rows vanish.  In the latter case, a row of $A$ must
vanish by Theorem \ref{AJTlargechar}.  In the former case, at
most $2$ entries from every row of $A$ can be nonzero.  If any
entry of a column is zero, either that column is identically $0$
or a row of $A$ iis identically $0$ by the same argument as in
Theorem \ref{AJTlargechar}.  Therefore, we can assume that exactly
two columns of $A$ are nonzero, and so we can reindex $A$ and 
normalize so that the first column is identically $1$.  Then
the coefficient $c_{i,j,0,\dotsc,0}$ is the $j^{\th}$ elementary
symmetric function of the entries in the second column.  Since
these are $0$ for every $0<j<d$ and 
$$c_{d,0,\dotsc,0}  + c_{0,d,\dotsc,0} = 0$$
so $c_{0,d,\dotsc,0} = -1$, the entries of the second column
of $A$ are the roots of the polynomial $z^{p-1} - 1$, which
are all the elements of the multiplicative group $k^{*}$.
\end{proof}

In fact, this result holds without the assumption that $q+1$ is
a prime.  This can be proven by using the proof above to show that
after normalizing one column, the other columns must contain all the
nonzero elements of the field.  Then, assuming that the first column
is identically $1$, and the next two columns are nonzero, computing
$c_{\alpha}$ for $\alpha = (i,q-i,1,0,\dotsc,0)$ for each $i$, one
can show via the Vandermonde determinant that either column $2$ or $3$
must be zero.

\section{Conclusion}
Proving Conjecture \ref{mainconj} it would be a great step forward in 
understanding permanental ideals and permanents in general.  For example, 
it would make it relatively easy to decide whether a matrix of size 
$d \x n$ where $n \geq 2d-1$ had vanishing maximal permanents.

It is somewhat striking that the only techniques used in this paper rely
on the multilinearity of permanents with some help from symmetry.  This
suggests that the multilinearity techniques could be successfully applied
to ideals defined by other multilinear functions on the $d \x d$ submatrices.
I believe that this could be a fruitful approach, but have not had time to
pursue this direction.

Another direction I have not pursued is the full primary decomposition of the
ideals.  Proposition \ref{34prop} (proved using Macaulay 2) says that $I_3(3,4)$
is radical, so there are no embedded components in that case.  For $I_3(3,5)$
on the other hand, there are over 200 embedded components, and the decomposition
is at this point inscrutable.

\providecommand{\bysame}{\leavevmode\hbox to3em{\hrulefill}\thinspace}
\providecommand{\MR}{\relax\ifhmode\unskip\space\fi MR }
\providecommand{\MRhref}[2]{%
  \href{http://www.ams.org/mathscinet-getitem?mr=#1}{#2}
}
\providecommand{\href}[2]{#2}


\begin{thebibliography}{1}

\bibitem{AT}
N.~Alon and M.~Tarsi, \emph{A nowhere-zero point in linear mappings},
  Combinatorica \textbf{9} (1989), no.~4, 393--395.

\bibitem{BBLS}
R.~D. Baker, J.~Bonin, F.~Lazebnik, and E.~Shustin, \emph{On the number of
  nowhere zero points in linear mappings}, Combinatorica \textbf{14} (1994),
  no.~2, 149--157. \MR{MR1289069 (95k:11160)}

\bibitem{M2}
Daniel~R. Grayson and Michael~E. Stillman, \emph{Macaulay 2, a software system
  for research in algebraic geometry}, Available at
  http://www.math.uiuc.edu/Macaulay2/.

\bibitem{GPS01}
G.-M. Greuel, G.~Pfister, and H.~Sch\"onemann, \emph{{\sc Singular} 2.0}, {A
  Computer Algebra System for Polynomial Computations}, Centre for Computer
  Algebra, University of Kaiserslautern, 2001, {\tt
  http://www.singular.uni-kl.de}.

\bibitem{LS}
Reinhard Laubenbacher and Irena Swanson, \emph{Permanental ideals}, J. Symbolic
  Computation \textbf{30} (2000), 195--295.

\bibitem{SPE}
Bernd Sturmfels, \emph{Solving systems of polynomial equations}, CBMS Regional
  Conference Series in Mathematics, vol.~97, Published for the Conference Board
  of the Mathematical Sciences, Washington, DC, 2002. \MR{2003i:13037}

\bibitem{Yu}
Yang Yu, \emph{The permanent rank of a matrix}, J. Combin. Theory Ser. A
  \textbf{85} (1999), no.~2, 237--242. \MR{MR1673948 (99j:15013)}

\end{thebibliography}
\end{document}